\documentclass[12pt]{article}
\usepackage{t1enc}
\usepackage{sw20elba}
\usepackage{amsmath}
\usepackage{amsfonts}
\usepackage{amssymb}
\usepackage{graphicx}
\newtheorem{theorem}{Theorem}

\newtheorem{definition}[theorem]{Definition}

\newtheorem{proposition}[theorem]{Proposition}

\newenvironment{proof}[1][Proof]{\textbf{#1.} }{\ \rule{0.5em}{0.5em}}

\begin{document}

\author{J\k{e}drzej \'{S}niatycki\\Department of Mathematics and Statistics\\University of Calgary\\Calgary, Alberta, Canada\\e-mail: sniat@math.ucalgary.ca}
\title{Singular Reduction for Nonlinear Control Systems}
\date{}
\maketitle
\begin{abstract}
We discuss smooth nonlinear control systems with symmetry. For a free and
proper action of the symmetry group, the reduction of symmetry gives rise to a
reduced smooth nonlinear control system. If the action of the symmetry group
is only proper, the reduced nonlinear control system need not be smooth.

Using the smooth calculus on nonsmooth spaces, provided by the theory of
differential spaces of Sikorski, we prove a generalization of Sussmann's
theorem on orbits of families of smooth vector fields.
\end{abstract}

\bigskip

\bigskip



\section{Introduction}

Geometric control theory is formulated in terms of smooth manifolds in order
to be able to use smooth calculus in the analysis of the system. The bundle
picture of nonlinear control was introduced by Brockett, \cite{Brockett},
because the local nonlinear control system dynamics on a manifold $M$, given
by
\[
\dot{x}=\varphi(x,u),
\]
where $\varphi:M\times U\rightarrow TM,$ was not an adequate description in
the case where the inputs depend on the states or the time histories of the states.

The role of symmetries was discussed by van der Schaft, \cite{Schaft}, and
Grizzle and Marcus, \cite{Grizzle-Marcus}, see also \cite{Bloch-Leonard} and
references quoted there. If the action of the symmetry group $G$ of the phase
space $M$ of the system is free and proper, and the control vector fields are
$G$-invariant, then the reduction of symmetry gives rise to a control system
on the space $\bar{M}=M/G$ of orbits of $G$ on $M$. The assumption that the
action of $G$ on $M$ is free and proper implies that $\bar{M}$ is a quotient
manifold of $M$. This is the setting for a discussion of symmetries in most of
the papers mentioned above.

If the action of $G$ on $M$ is not free, then the orbit space $\bar{M}=M/G$
may have singularities. In the case of a proper action the orbit space is a
stratified space, \cite{DK}. The reduction of a proper action of a symmetry
group of a Hamiltonian system has been a subject of investigation for several
years. A convenient tool to study reduced spaces, which are not manifolds, is
the theory of differential spaces of Sikorski, \cite{sikorski1}. It was first
used implicitly by Sjamaar and Lerman, \cite{SL}. A differential space
formulation of singular reduction is given in \cite{CS} for unconstrained
Hamiltonian systems, and in \cite{Sniatycki2001} for nonholonomically
constrained Hamiltonian system.

The aim of this paper is to discuss reduction of symmetries of nonlinear
control systems, and to show how the calculus on differential spaces enables
us to get results in a non-smooth setting. For a free and proper action of the
symmetry group, the reduction reproduces the bundle formulation. If the action
of the symmetry group is only proper, we get a differential space analogue of
the bundle formulation of Brockett. Nevertheless, we are able to prove a
generalization of Sussmann's Theorem on orbits of families of local vector
fields, \cite{Sussmann}.

In Section 2, we describe the bundle formulation for smooth nonlinear control
systems, and their symmetries. Section 3 is devoted to a discussion of
symmetry reduction. A generalization of Sussmann's Theorem is proved in
Section 4. Section 5 contains concluding remarks. A review of techniques of
differential space theory is given in Section 6

\section{Nonlinear Control Systems with Symmetry}

A smooth nonlinear control system is defined to be a quadruple $(B,M,\pi
,\varphi)$ such that

\begin{description}
\item [(i)]$(B,M,\pi)$ is a fibre bundle with total space $B$, base space $M$
and projection $\pi:B\rightarrow M,$ and

\item[(ii)] $\varphi:B\rightarrow TM$ is a bundle morphism such that, for each
$x\in M$ and each $b\in B_{x}=\pi^{-1}(x)$, $\varphi(b)\in T_{x}M$.
\end{description}

The assumption that $(B,M,\pi)$ is a fibre bundle implies that there exists a
family $\Gamma(M,B)$ of smooth local sections $\sigma$ of $\pi:B\rightarrow M$
such that $M$ is covered by the domains of $\sigma\in\Gamma(M,B).$ For each
$\sigma\in\Gamma(M,B)$, the composition $X=\varphi\raisebox{2pt}%
{$\scriptstyle\circ\, $}\sigma$ is a control vector field on $M$. In this way
we obtain a family $\mathcal{D}$
\begin{equation}
\mathcal{D}=\{X=\varphi\raisebox{2pt}{$\scriptstyle\circ\, $}\sigma\mid
\sigma\in\Gamma(M,B)\}\label{F}%
\end{equation}
of locally defined vector fields on $M$ such that $M$ is covered by the
domains of $X\in\mathcal{D}.$ By choosing $\Gamma(M,B)$ we have effectively
reduced the non-linear control system $(B,M,\pi,\varphi)$ to a piecewise
linear system given by the family $\mathcal{D}$ of local vector fields on $M$.

An example of a control problem on $M$ is the analysis of the structure of
accessible sets of a family $\mathcal{D}$. For each $X\in\mathcal{D}$, we
denote by $\exp(tX)$ the local one-parameter local group of diffeomorphisms of
$M$ generated by $X$. For every $x\in M$, the accessible set of $\mathcal{D}$
through $x$ is
\begin{equation}
N_{x}=\{\exp(t_{n}X_{n})\raisebox{2pt}{$\scriptstyle\circ\, $}...\raisebox
{2pt}{$\scriptstyle\circ\, $}\exp(t_{1}X_{1})\mid n\in\mathbb{N},\text{ }%
t_{1},...,t_{n}\in\mathbb{R},\text{ }X_{1},...,X_{n}\in\mathcal{D}%
\}.\label{accessible}%
\end{equation}
It has been shown by Sussmann that $N_{x}$ is a manifold immersed in $M,$
\cite{Sussmann}. The family of accessible sets of $\mathcal{D}$ defines on $M
$ the structure of a smooth foliation with singularities, \cite{Stefan}.

Let $G$ be a Lie group with Lie algebra $\frak{g}$, and let
\[
\Theta:G\times B\rightarrow B:(g,u)\mapsto\Theta(g,u)\equiv\Theta_{g}(u)\equiv
gu,
\]
and
\[
\Phi:G\times M\rightarrow M:(g,x)\mapsto\Phi(g,x)\equiv\Phi_{g}(u)\equiv gx,
\]
be left actions of $G$ on $B$, and $M$, respectively. We say that $G$ is a
symmetry group of the control system $(B,M,\pi,f)$ if the map $\varphi
:B\rightarrow TM$ intertwines the action $\Theta$ on $B$ and $T\Phi$ on $TM $.
In other words, $G$ is a symmetry if, for every $g\in G$,
\[
\varphi\raisebox{2pt}{$\scriptstyle\circ\, $}\Theta_{g}=T\Phi_{g}%
\raisebox{2pt}{$\scriptstyle\circ\, $}\varphi.
\]

We consider here a special case in which all sections $\sigma\in\Gamma(M,B)$
intertwine the action $\Phi$ on $M$ and the action $\Theta$ on $B$. In other
words, we assume that
\begin{equation}
\sigma\raisebox{2pt}{$\scriptstyle\circ\, $}\Phi_{g}=\Theta_{g}\raisebox
{2pt}{$\scriptstyle\circ\, $}\sigma\label{invariance}%
\end{equation}
for every $\sigma\in\Gamma(M,B)$ and $g\in G$. In this case all vector fields
$X\in\mathcal{F}$ are $G$-invariant. Thus, we are dealing with a control
system on a manifold with symmetry in which all controls have the same symmetry.

\section{Reduction}

In this section we describe the control system obtained by symmetry reduction
of a smooth nonlinear control system with symmetry.

If the actions $\Theta$ and $\Phi$ are free and proper, then orbit spaces
$\bar{B}=B/G$ and $\bar{M}=M/G$ are quotient manifolds of $B$ and $M$,
respectively, with projection maps $\beta:B\rightarrow\bar{B}$ and
$\mu:M\rightarrow\bar{M}.$ Since the map $\pi:B\rightarrow M$ intertwines the
actions $\Theta$ and $\Phi$, it induces a map $\bar{\pi}:\bar{B}%
\rightarrow\bar{M}$ such that
\[
\mu\raisebox{2pt}{$\scriptstyle\circ\, $}\pi=\bar{\pi}\raisebox{2pt}%
{$\scriptstyle\circ\, $}\beta.
\]

Let $\psi=T\mu\raisebox{2pt}{$\scriptstyle\circ\, $}\varphi:B\rightarrow
T\bar{M}$. For every $g\in G$, and $u\in B$,
\[
\psi(gu)=T\mu(\varphi(\Theta_{g}u))=T\mu(T\Phi_{g}(\varphi(u))=T\mu
(\varphi(u))=\psi(u).
\]
Thus, $\psi$ is constant on orbits of $G$, and it pushes forward to a smooth
map $\bar{\varphi}:\bar{B}\rightarrow T\bar{M}$ such that
\[
T\mu\raisebox{2pt}{$\scriptstyle\circ\, $}\varphi=\bar{\varphi}\raisebox
{2pt}{$\scriptstyle\circ\, $}\beta.
\]

\begin{proposition}
The quadruple $(\bar{B},\bar{M},\bar{\pi},\bar{\varphi})$ is a smooth
nonlinear control system.
\end{proposition}

\begin{proof}
Since $\bar{\pi}\raisebox{2pt}{$\scriptstyle\circ\, $}\beta=\mu\raisebox
{2pt}{$\scriptstyle\circ\, $}\pi,$ it follows that $\beta^{-1}(\bar{\pi}%
^{-1}(V))=\pi^{-1}(\mu^{-1}(V))$ for every $V\subseteq\bar{M}.$ Hence,
$\bar{\pi}^{-1}(V)=\beta(\pi^{-1}(\mu^{-1}(V)).$

Since the actions $\Theta$ and $\Phi$ are free and proper, they introduce
structure of\ a (left) $G$-principal fibre bundle on $\beta:B\rightarrow
\bar{B}$ and $\mu:M\rightarrow\bar{M}$, respectively. Hence, for every
$\bar{x}\in\bar{M}$, there exists a neigbourhood $V$ of $\bar{x}$ in $\bar{M}
$, such that $\mu^{-1}(V)\cong G\times V$. Let $e$ denote the identity in $G
$. Since $\pi:B\rightarrow M$ is locally trivial, there exists a neighbourhood
of $x=(e,\bar{x})\in\mu^{-1}(V)\subseteq M$ of the form $U\times S_{x}$, where
$U$ is an open neigbourhood of $e$ in $G$ and $S_{x}$ is a slice through $x$
for the action of $G$ on $M,$ such that
\[
\pi^{-1}(U\times S_{x})\cong\pi^{-1}(x)\times U\times S_{x}.
\]
By shrinking $V$ and $S_{x}$, if necessary, we may assume that $V=\mu(S_{x}).$
This implies that
\[
\pi^{-1}(\mu^{-1}(V))\cong\pi^{-1}(G\times V)\cong\pi^{-1}(x)\times G\times
S_{x}.
\]
Hence,
\[
\bar{\pi}^{-1}(V)\cong\beta(\pi^{-1}(\mu^{-1}(V))\cong\beta(\pi^{-1}(x)\times
G\times S_{x})=\pi^{-1}(x)\times S_{x}\cong\pi^{-1}(x)\times V\text{.}
\]
Moreover,
\begin{align*}
\bar{\pi}^{-1}(\bar{x})  & =\{\bar{b}\in\bar{B}\mid\bar{\pi}(\bar{b})=\bar
{x}\}\cong\{Gb\mid b\in B,\text{ }\beta\raisebox{2pt}{$\scriptstyle\circ\,
$}\pi(b)=\bar{x}\}\\
& \cong\{(b^{\prime},Gg^{\prime},x)\subset\pi^{-1}(x)\times G\times S_{x}\mid
g^{\prime}\in G\}\\
& \cong\{(b^{\prime},x)\in\pi^{-1}(x)\times\{x\}\}\cong\pi^{-1}(x).
\end{align*}
This implies that $\bar{\pi}:\bar{B}\rightarrow\bar{M}$ is locally trivial.

The map $\bar{\varphi}:\bar{B}\rightarrow T\bar{M}$ satisfies $T\mu
\raisebox{2pt}{$\scriptstyle\circ\, $}\varphi=\bar{\varphi}\raisebox
{2pt}{$\scriptstyle\circ\, $}\beta.$ Hence, for each $\bar{x}\in\bar{M},$
$\bar{b}\in\bar{B}_{\bar{x}}=\bar{\pi}^{-1}(\bar{x})$, and $b\in\beta
^{-1}(\bar{b})$, we have
\[
\bar{\varphi}(\bar{b})=\bar{\varphi}(\beta(b))=T\mu(\varphi(b))\in T\mu
(T_{\pi(b)}M)\subset T_{\mu(\pi(b))}\bar{M}=T_{\bar{x}}\bar{M}
\]
because $\bar{x}=\bar{\pi}(\bar{b})=\bar{\pi}(\mu(b))=\mu(\pi(b))$. This
completes the proof.
\end{proof}

Let $\Gamma(\bar{M},\bar{B})$ denote a space of smooth sections $\bar{\sigma
}:\bar{M}\rightarrow\bar{B}$ of $\bar{\pi}:\bar{B}\rightarrow\bar{M}$. For
each $\bar{\sigma}\in\Gamma(\bar{M},\bar{B})$, the composition $\bar{X}%
=\bar{f}\raisebox{2pt}{$\scriptstyle\circ\, $}\bar{\sigma}$ is a control
vector field on $\bar{M}$. As before, we obtain a family $\overline
{\mathcal{D}}$ of vector fields on $\bar{M}$ parametrized by $\bar{\sigma}%
\in\Gamma(\bar{M},\bar{B}).$ In other words,
\begin{equation}
\overline{\mathcal{D}}=\{\bar{X}=\bar{\varphi}\raisebox{2pt}{$\scriptstyle
\circ
\, $}\bar{\sigma}\mid\bar{\sigma}\in\Gamma(\bar{M},\bar{B})\}.\label{Fbar}%
\end{equation}

\begin{proposition}
For a free and proper action of $G$ on $M$, given a vector field $\bar{X}$ on
the orbit space $\bar{M}=M/G$ there exists a $G$-invariant vector field $X $
on $M$ such that $\bar{X}=\mu_{\ast}X$, where $\mu:M\rightarrow\bar{M}$ is the
orbit map.
\end{proposition}

\begin{proof}
Let $\mathrm{hor}TM$ be a connection on the principal $G$-bundle
$\mu:M\rightarrow\bar{M}$. That is, $\mathrm{hor}TM$ is a $G$-invariant
distribution on $M$ such that
\[
\mathrm{hor}TM\oplus\ker T\mu=TM.
\]
Let $X$ be the horizontal lift of $\bar{X}$. In other words, $X$ is a vector
field on $M$ with values in $\mathrm{hor}TM$ such that $\bar{X}=\mu_{\ast}X $.
Since $\mathrm{hor}TM$ is $G$-invariant, it follows that $X$ is $G$-invariant.
\end{proof}

It follows that the piece-wise linearization of the reduced smooth non-linear
control system $(\bar{B},\bar{M},\bar{\pi},\bar{\varphi})$ corresponds to a
piece-wise linearization of the original smooth non-linear control system
$(B,M,\pi,\varphi)$ in terms of a family $\mathcal{D}$ of $G $-invariant
vector fields on $M$.

If the actions $\Theta$ and $\Phi$ are not free, the orbit spaces $\bar{B}$
and $\bar{M}$ need not be manifolds. If $\Theta$ and $\Phi$ are proper, then
$\bar{B}$ and $\bar{M}$ are stratified spaces, \cite{DK}. Following Schwarz,
\cite{Schwarz}, we define differential structures on $\bar{B}$ and $\bar{M}$
in terms of $G$-invariant smooth functions on $B$ and $M$, respectively. More
precisely,
\[
C^{\infty}(\bar{B})=\{h:\bar{B}\rightarrow\mathbb{R}\mid h\raisebox
{2pt}{$\scriptstyle\circ\, $}\beta\in C^{\infty}(B)\},
\]
and
\[
C^{\infty}(\bar{M})=\{h:\bar{M}\rightarrow\mathbb{R}\mid h\raisebox
{2pt}{$\scriptstyle\circ\, $}\mu\in C^{\infty}(M)\}.
\]
The spaces $\bar{B}$ and $\bar{M}$ endowed with these differential structures
are Hausdorff differential spaces in the sense of Sikorski, \cite{Sikorski}.
In \cite{CS}, it has been shown that orbit spaces of a proper action are
differential spaces that are locally diffeomorphic to subsets of
$\mathbb{R}^{n}.$ Such spaces have been introduced by Aronszajn,
\cite{Aronszajn}, under the name of subcartesian spaces. In the appendix, we
review properties of smooth subcartesian spaces following \cite{Sniatycki}.

As in the case of a free and proper action, we have a smooth projection
$\bar{\beta}:\bar{B}\rightarrow\bar{M}.$ In order to describe the mapping
$\bar{\varphi}:\bar{B}\rightarrow T\bar{M}$, we have to define what we mean
here by the ``tangent bundle space'' of a subcartesian space. This problem has
been a subject of many papers, see \cite{Andrzejczak},
\cite{aronszajn-szeptycki}, \cite{Kowalczyk}, \cite{marshall},
\cite{marshall1}, \cite{Matuszczyk}, \cite{Matuszczyk1},
\cite{Matuszczyk-Waliszewski}, \cite{sikorski2}, \cite{spallek2},
\cite{Wierzbicki} and references quoted there. Different notions of tangent
vectors, which are equivalent on a manifold, are need not be equivalent in the
case of a differential space.

We begin with the Zariski tangent bundle space $T^{Z}M$. For $x\in M$, the
Zariski tangent space $T_{x}^{Z}\bar{M}$ consists derivations at $x$ of
$C^{\infty}(\bar{M}).$ In other words, an element of $T_{x}^{Z}\bar{M}$ is a
linear map
\[
v:C^{\infty}(\bar{M})\rightarrow\mathbb{R}:h\mapsto v\cdot h
\]
satisfying Leibniz' rule
\[
v\cdot(h_{1}h_{2})=(v\cdot h_{1})h_{2}(x)+h_{1}(x)(v\cdot h_{2})
\]
for all $h_{1},h_{2}\in C^{\infty}(\bar{M}).$

Vectors in $TM$ act by derivation on $C^{\infty}(M)$. If $u\in B,$ and
$x=\pi(u)$, then $\varphi(u)\in T_{x}M$ acts on $f\in C^{\infty}(M)$ by
\[
\varphi(u)\cdot f=\frac{d}{dt}f(c(t))_{\mid t=0},
\]
where $t\mapsto c(t)$ is a curve in $M$ such that $c(0)=x$ and $\dot
{c}(0)=\varphi(u).$ For every $g\in G$,
\[
\varphi(\Theta_{g}u)\cdot f=T\Phi_{g}(\varphi(u))\cdot f=\varphi\cdot\Phi
_{g}^{\ast}f=\frac{d}{dt}f(\Phi_{g}(c(t)))_{\mid t=0}.
\]
If $f$ is $G$-invariant, then $f\raisebox{2pt}{$\scriptstyle\circ\, $}\Phi
_{g}=f$, and $\varphi(\Theta_{g}u)\cdot f=\varphi(u)\cdot f$ for every $g\in
G$. In this case, $\varphi(u)\cdot f$ depends only on $\bar{u}=\beta(u)\in
\bar{B}.$ Since every $G$-invariant function on $M$ is of the form
$f=h\raisebox{2pt}{$\scriptstyle\circ\, $}\mu,$ for a unique $h\in C^{\infty
}(\bar{M})$, we have a map $\bar{\varphi}:\bar{B}\rightarrow T^{Z}\bar{M}$
such that
\begin{equation}
\bar{\varphi}(\bar{u})\cdot h=\beta(u)\cdot(h\raisebox{2pt}{$\scriptstyle
\circ\, $}\mu)\label{phibar}%
\end{equation}
for every $h\in C^{\infty}(\bar{M})$, where $u$ is any element of $\beta
^{-1}(\bar{u})$.

Before we can claim that $(\bar{B},\bar{M},\bar{\pi},\bar{\varphi})$ is the
reduced control system, we have to examine the role played by $\varphi$ and
$\bar{\varphi}$. In the preceding section, we have used the map $\varphi
:B\rightarrow TM$ to associate to a family $\Gamma(M,B)$ of local sections of
$\pi:B\rightarrow M$ a family $\mathcal{D}$ of locally defined vector fields
on $M$. The definition of accessible sets of $\mathcal{D}$ is based on the
fact that every vector field $X$ gives rise to a local one-parameter group of
local diffeomorphisms $\exp(tX).$ By assumption, every section $\sigma
\in\Gamma(M,B)$ intertwines the actions of $G$ on $M$ and $B$, see equation
(\ref{invariance}). Hence, it gives rise to a local section $\bar{\sigma}$ of
$\bar{\pi}:\bar{B}\rightarrow\bar{M}$ such that
\[
\bar{\sigma}\raisebox{2pt}{$\scriptstyle\circ\, $}\beta=\mu\raisebox
{2pt}{$\scriptstyle\circ\, $}\sigma.
\]
Moreover, the local vector field $X=\varphi\raisebox{2pt}{$\scriptstyle
\circ\, $}\sigma$ is $G$-invariant. Hence, it gives rise to a local derivation
$\bar{X}$ of $C^{\infty}(\bar{M})$ such that, for every $x\in\mathrm{domain}X$
and $h\in C^{\infty}(\bar{M})$
\[
\bar{X}(\mu(x))\cdot h=X(x)\cdot(h\raisebox{2pt}{$\scriptstyle\circ\, $}\mu).
\]
In this way we obtain a family $\overline{\mathcal{D}}$ of local derivations
of $C^{\infty}(\bar{M})$. Moreover, for every $\bar{X}\in\overline
{\mathcal{D}}$ and $\bar{x}\in\mathrm{domain}\bar{X}$,
\[
\bar{X}(x)=\bar{\varphi}(\bar{\sigma}(x)).
\]

Hence, the assumption (\ref{invariance}) implies that the family $\Gamma(M,B)$
of local sections of $\pi:B\rightarrow M$ gives rise to a family $\Gamma
(\bar{M},\bar{B})$ of local sections of $\bar{\pi}:\bar{B}\rightarrow\bar{M}$.
Moreover, the family $\overline{\mathcal{D}}$ of local derivations of $\bar
{M}$ is given by
\[
\overline{\mathcal{D}}=\{\bar{X}=\bar{\varphi}\raisebox{2pt}{$\scriptstyle
\circ
\, $}\bar{\sigma}\mid\bar{\sigma}\in\Gamma(\bar{M},\bar{B})\}.
\]

The question arises if local derivations $\bar{X}\in\overline{\mathcal{D}}$
generate local one-parameter groups of local diffeomorphisms of $\bar{M}$.
Every local derivation of $C^{\infty}(\bar{M})$ extends locally to a global
derivation of $C^{\infty}(\bar{M}),$ \cite{Sniatycki}. However, not all global
derivations of $C^{\infty}(\bar{M})$ generate local one-parameter groups of
local diffeomorphisms of $\bar{M}$. Global derivations of $C^{\infty}(\bar
{M})$ that generate local one-parameter groups of local diffeomorphisms of
$\bar{M}$ are called in \cite{Sniatycki} vector fields on $\bar{M},$ for
details see Appendix. In the remainder of this section we show that all local
derivations $\bar{X}\in\overline{\mathcal{D}}$ are local vector fields on $M$
in the sense that they generate local one-parameter groups of local
diffeomorphisms of $\bar{M}$.

It follows from the existence of $G$-invariant partitions of unity for proper
actions, see \cite{BC}, that every locally defined $G$-invariant vector field
$X$ can be locally extended to a globally defined $G$-invariant vector field.
In other words, for every $x\in\mathrm{domain}X,$ there exists a neighbourhood
$U$ of $x\in M$ and a globally defined $G$-invariant vector field $X^{\prime}$
on $M$ such that the restrictions to $U$ of $X$ and $X^{\prime}$ coincide,
i.e. $X_{\mid U}=X_{\mid U}^{\prime}$. Hence, vector fields in $\mathcal{D}$
are locally restrictions of globally defined $G$-invariant vector fields on
$M$. Therefore, the local one-parameter group $\exp(tX)$ of local
diffeomorphisms of $\mathrm{domain}X$ and the local one-parameter group
$\exp(tX^{\prime})$ of local diffeomorphisms of $M$ coincide on $U.$

By construction, a vector field $\bar{X}\in\overline{\mathcal{D}}$ is the
push-forward by $\mu$ of a vector field $X\in\mathcal{D}.$ Since \ $X$ is
$G$-invariant, the local one-parameter group $\exp(tX)$ of local
diffeomorphisms of $\mathrm{domain}X\subseteq M$ generated by $X$ preserves
$G$-orbits in $\mathrm{domain}X$. Hence, it induces a local one-parameter
group $\exp(t\bar{X})$ of local transformations of the orbit space
$(\mathrm{domain}X)/G=\mathrm{domain}\bar{X}.$ For every $G$-invariant
function $f\in C^{\infty}(\mathrm{domain}X)$, the pull-back $\exp(tX)^{\ast
}f=f\raisebox{2pt}{$\scriptstyle\circ\,
$}\exp(tX)$ is $G$-invariant. Hence, $\exp(t\bar{X})^{\ast}h$ is smooth for
every $h\in C^{\infty}(\mathrm{domain}\bar{X})$. This implies that $\exp
(t\bar{X})$ is a local one-parameter group of local diffeomorphisms of
$\mathrm{domain}\bar{X}$.

On the other hand, $\bar{X}^{\prime}=\mu_{\ast}X^{\prime}$ is a derivation of
$C^{\infty}(\bar{M})$ and it induces a local one-parameter group $\exp
(t\bar{X}^{\prime})$ of local diffeomorphisms of $\bar{M}$ which coincides
with $\exp(t\bar{X})$ on $\bar{U}=U/G$. Moreover, for every point $x\in M$
there is a neighbourhood $U$ of $x$ in $M$, and a $G$-invariant vector field
$X^{\prime}$ on $M$ such that the above condition is satisfied. Hence, for
every $\bar{X}\in\mathcal{D}$, the local one-parameter group $\exp(t\bar{X})$
of local diffeomorphisms of $\mathrm{domain}\bar{X}$ is given locally by
restrictions of local one-parameter groups of local diffeomorphisms of $M$.

We denote by $T\bar{M}$ the set of values of all vector fields on $\bar{M}.$
In other words, for every $x\in\bar{M}$ and $\bar{v}\in T_{x}\bar{M}$, there
exists a global derivation $\bar{X}$ of $C^{\infty}(\bar{M})$, generating a
local one-parameter group $\exp(t\bar{X})$ of local diffeomorphisms of
$\bar{M}$ such that $\bar{X}(x)=\bar{v}.$ The discussion above implies that
the map $\bar{\varphi}$ defined by equation (\ref{phibar}) has values in
$T\bar{M}.$ \ We shall see in the next section that this property implies that
orbit spaces of $\mathcal{\bar{D}}$ are smooth manifolds. In the following, we
write $\bar{\varphi}:\bar{B}\rightarrow T\bar{M}$, and refer to the quadruple
$(\bar{B},\bar{M},\bar{\pi},\bar{\varphi})$ as the reduced control system.

\section{Generalized Sussmann's Theorem}

We prove here a generalization of Sussmann's Theorem, \cite{Sussmann}, to
subcartesian spaces.

Let $\bar{M}$ be a subcartesian space and \thinspace$\mathcal{X}(\bar{M})$ the
family of all globally defined vector fields on $\bar{M}$. In other words,
$X\in\mathcal{X}(\bar{M})$ if and only if $X$ is a derivation of $C^{\infty
}(\bar{M})$ that generates a local one-parameter local group of
diffeomorphisms of $\bar{M}.$

\begin{theorem}
Let $\overline{\mathcal{D}}$ be a family of local vector fields on a
subcartesian space $\bar{M}$ such that, for every $\bar{X}\in\overline
{\mathcal{D}}$ and $\bar{x}\in\mathrm{domain}\bar{X},$ there exists a
neighbourhood $\bar{U}$ of $\bar{x}$ and a vector field $\bar{X}^{\prime}%
\in\mathcal{X}(\bar{M})$ such that $\bar{X}_{\mid\bar{U}}=\bar{X}_{\mid\bar
{U}}^{\prime}.$ For each $\bar{x}\in\bar{M}$ the orbit $\bar{N}_{\bar{x}}$ of
$\overline{\mathcal{D}}$ through $\bar{x},$ defined by
\[
\bar{N}_{\bar{x}}=\{\exp(t_{n}\bar{X}_{n})\raisebox{2pt}{$\scriptstyle\circ
\, $}...\raisebox{2pt}{$\scriptstyle\circ\, $}\exp(t_{1}\bar{X}_{1})\mid
n\in\mathbb{N},\text{ }t_{1},...,t_{n}\in\mathbb{R},\text{ }\bar{X}%
_{1},...,\bar{X}_{n}\in\overline{\mathcal{D}}\}
\]
is a smooth manifold.
\end{theorem}

\begin{proof}
For each $\bar{X}_{i}\in\overline{\mathcal{D}}$ in the expression for $\bar
{N}_{\bar{x}}$, we can replace $\exp(t_{i}\bar{X}_{i})$ by $\exp(t_{i}\bar
{X}_{i}^{\prime})$, where $\bar{X}_{i}^{\prime}\in\mathcal{X}(\bar{M})$. This
implies that $\bar{N}_{\bar{x}}$ is contained in the orbit
\[
\bar{M}_{\bar{x}}=\{\exp(t_{n}\bar{X}_{n}^{\prime})\raisebox{2pt}%
{$\scriptstyle\circ\, $}...\raisebox{2pt}{$\scriptstyle\circ\, $}\exp
(t_{1}\bar{X}_{1}^{\prime})\mid n\in\mathbb{N},\text{ }t_{1},...,t_{n}%
\in\mathbb{R},\text{ }\bar{X}_{1}^{\prime},...,\bar{X}_{n}^{\prime}%
\in\mathcal{X}(\bar{M})\}
\]
of $\mathcal{X}(\bar{M})$ passing through $x$. It has been proved in
\cite{Sniatycki} that $\bar{M}_{\bar{x}}$ is a smooth manifold.

Let $\overline{\mathcal{D}}_{\mid\bar{M}_{x}}$ denote the family of local
vector fields on $\bar{M}_{\bar{x}}$ obtained by the restriction to $\bar
{M}_{\bar{x}}$ of local vector fields in $\overline{\mathcal{D}}$. We can
write
\[
\bar{N}_{\bar{x}}=\{\exp(t_{n}\bar{X}_{n})\raisebox{2pt}{$\scriptstyle\circ
\, $}...\raisebox{2pt}{$\scriptstyle\circ\, $}\exp(t_{1}\bar{X}_{1})\mid
n\in\mathbb{N},\text{ }t_{1},...,t_{n}\in\mathbb{R},\text{ }\bar{X}%
_{1},...,\bar{X}_{n}\in\overline{\mathcal{D}}_{\mid\bar{M}_{\bar{x}}}\}.
\]
By Sussmann's Theorem, \cite{Sussmann}, $\bar{N}_{\bar{x}}$ is an immersed
submanifold of $\bar{M}_{\bar{x}}.$ Hence, $\bar{N}_{\bar{x}}$ is a smooth manifold.
\end{proof}

\section{Concluding remarks}

It follows from the discussion in the preceding sections that the notion of a
smooth nonlinear control system $(B,M,\pi,\varphi)$ can be naturally extended
to the case when the spaces $B$ and $M$ have singularities. Reduction of
symmetries gives rise to stratified spaces with relatively mind singularities.
However, as illustrated in the proof of Theorem 3, most of arguments used here
are valid for subcartesian spaces. Since subcartesian spaces are locally
diffeomorphic to arbitrary subsets of $\mathbb{R}^{n}$, it follows that
allowable singularities are restricted only by finiteness of dimension.

\section{Appendix: Differential spaces}

Differential spaces were introduced by Sikorski, \cite{sikorski1}, see also
\cite{sikorski2} and \cite{Sikorski}. Their structure has been investigated by
several authors, see \cite{Andrzejczak}, \cite{aronszajn-szeptycki},
\cite{Kowalczyk}, \cite{marshall}, \cite{marshall1}, \cite{Matuszczyk},
\cite{Matuszczyk1}, \cite{Matuszczyk-Waliszewski}, \cite{spallek1},
\cite{spallek2}, \cite{Walczak}, \cite{Wierzbicki} and references quoted there.

A differential structure on a topological space $S$ is a family of functions
$C^{\infty}(S)$ satisfying the following conditions:

\begin{description}
\item [2.1.]The family
\[
\{f^{-1}((a,b))\mid f\in C^{\infty}(S),\,a,b\in\mathbb{R}\}
\]
is a sub-basis for the topology of $S.$

\item[2.2.] If $f_{1},...,f_{n}\in C^{\infty}(S)$ and $F\in C^{\infty
}(\mathbb{R}^{n})$, then $F(f_{1},...,f_{n})\in C^{\infty}(S).$

\item[2.3.] If $f:S\rightarrow\mathbb{R}$ is such that, for every $x\in S$,
there exist an open neighbourhood $U_{x}$ of $x$ and a function $f_{x}\in
C^{\infty}(S)$ satisfying
\[
f_{x}\mid U_{x}=f\mid U_{x},
\]
then $f\in C^{\infty}(S).$ Here the vertical bar $\mid$ denotes the restriction.
\end{description}

A differential space is a topological space endowed with a differential
structure.\ Let $R$ and\thinspace$S$ are differential spaces with differential
structures $C^{\infty}(R)$ and $C^{\infty}(S),$ respectively. A map
$\rho:R\rightarrow S$ is said to be smooth if $\rho^{\ast}f\in C^{\infty}(R)$
for all $f\in C^{\infty}(S)$. A smooth map between differential spaces is a
diffeomorphism if it is invertible and its inverse is smooth.\ 

Clearly, smooth manifolds are differential spaces. However, the category of
differential spaces is much larger than the category of manifolds.

If $R$ is a differential space with differential structure $C^{\infty}(R)$ and
$S$ is a subset of $R$, then we can define a differential structure
$C^{\infty}(S)$ on $S$ as follows. A function $f:S\rightarrow\mathbb{R}$ is in
$C^{\infty}(S)$ if and only if, for every $x\in S,$ there is an open
neighborhood $U$ of $x$ in $R$ and a function $f_{x}\in C^{\infty}(R)$ such
that $f|(S\cap U)=f_{x}|(S\cap U)$. The differential structure $C^{\infty}(S)$
described above is the smallest differential structure on $S$ such that the
inclusion map $\iota:S\rightarrow R$ is smooth. We shall refer to $S$ with the
differential structure $C^{\infty}(S)$ described above as a differential
subspace of $R$. If $S$ is a closed subset of $R,$ then the differential
structure $C^{\infty}(S)$ described above consists of restrictions to $S$ of
functions in $C^{\infty}(R)$.

A differential space $R$ is said to be locally diffeomorphic to a differential
space $S$ if, for every $x\in R$, there exists a neighbourhood $U$ of $x$
diffeomorphic to an open subset $V$ of $S$. More precisely, we require that
the differential subspace $U$ of $R$ be diffeomorphic to the differential
subspace $V$ of $S$. A differential space $R$ is a smooth manifold of
dimension $n$ if and only if it is locally diffeomorphic to $\mathbb{R}^{n}$.

A Hausdorff differential space that is locally diffeomorphic to a subset of
$\mathbb{R}^{n}$ is called a subcartesian space. The original definition of
subcartesian space was given by Aronszajn in terms of a singular atlas,
\cite{Aronszajn}, see also\cite{aronszajn-szeptycki1},
\cite{aronszajn-szeptycki}, \cite{marshall1} and \cite{marshall}. The
characterization of subcartesian spaces used here can be found in
\cite{spallek1} and \cite{Walczak}. In the following we review properties of
families of vector fields on subcartesian spaces. Proofs of theorems stated
here can be found in \cite{Sniatycki}.

Let $S$ be a subcartesian space with a differential structure $C^{\infty}(S)
$. A derivation on $C^{\infty}(S)$ is a linear map $X:C^{\infty}(S)\rightarrow
C^{\infty}(S):f\mapsto X\cdot f$ satisfying Leibniz' rule
\begin{equation}
X\cdot(f_{1}f_{2})=(X\cdot f_{1})f_{2}+f_{1}(X\cdot f_{2}).\label{Leibniz}%
\end{equation}
We denote the space of derivations of $C^{\infty}(S)$ by $\mathrm{Der}%
C^{\infty}(S).$ It has the structure of a Lie algebra with the Lie bracket
$[X_{1},X_{2}]$ defined by
\[
\lbrack X_{1},X_{2}]\cdot f=X_{1}\cdot(X_{2}\cdot f)-X_{2}\cdot(X_{1}\cdot f)
\]
for every $X_{1},X_{2}\in\mathrm{Der}C^{\infty}(S)$ and $f\in C^{\infty}(S).$

A local diffeomorphism $\varphi$ of $S$ to itself is a diffeomorphism
$\varphi:U\rightarrow V$, where $U$ and $V$ are open differential subspaces of
$S.$ For each $f\in C^{\infty}(S),$ the restriction of $f$ to $V$ is in
$C^{\infty}(V)$, and $\varphi^{\ast}f=f\raisebox{2pt}{$\scriptstyle\circ
\, $}\varphi$ is in $C^{\infty}(U)$. If $\varphi^{\ast}f$ coincides with the
restriction of $f$ to $U,$ we say that $f$ is $\varphi$-invariant, and write
$\varphi^{\ast}f=f.$ For each $X\in\mathrm{Der}(C^{\infty}(S)),$ the
restriction of $X$ to $U$ is in $\mathrm{Der}(C^{\infty}(U))$, and the
push-forward $\varphi_{\ast}X$ of $X$ by $\varphi$ is a derivation of
$C^{\infty}(V)$ such that
\begin{equation}
(\varphi_{\ast}X)\cdot(f\mid V)=\varphi^{-1\ast}(X\cdot(\varphi^{\ast
}f))\text{ for all }f\in C^{\infty}(S).\label{push-forward}%
\end{equation}
Since all functions in $C^{\infty}(V)$ locally coincide with restrictions to
$V$ of functions in $C^{\infty}(S)$, equation (\ref{push-forward}) determines
$\varphi_{\ast}X$ uniquely. If $\varphi_{\ast}X$ coincides with the
restriction of $X$ to $V$, we say that $X$ is $\varphi$-invariant and write
$\varphi_{\ast}X=X.$

Let $I$ be an interval in $\mathbb{R}$. A smooth map $c:I\rightarrow
S:t\mapsto x(t)$ is an integral curve of a derivation $X$ if
\[
\frac{d}{dt}f(x(t))=(X\cdot f)(x(t))
\]
for all $f\in C^{\infty}(S)$ and $t\in I$.

\begin{theorem}
For every derivation $X$ on a subcartesian space $S$ and each point $x\in S$,
there exists a unique maximal integral curve $c$ of $X$ such that $c(0)=x.$
\end{theorem}

\begin{definition}
A \emph{vector field} on a subcartesian space $S$ is a derivation $X$ of
$C^{\infty}(S)$ such that translations along integral curves of $X$ give rise
to local diffeomorphisms of $S.$
\end{definition}%

%TCIMACRO{\TeXButton{noindent}{\noindent}}%
%BeginExpansion
\noindent
%EndExpansion
There is a simple criterion characterizing vector fields on a subcartesian
space; namely,

\begin{theorem}
A derivation $X$ of $C^{\infty}(S)$ is a vector field on $S$ if and only if
the domains of maximal integral curves of $X$ are open in $\mathbb{R}.$
\end{theorem}

If $X$ is a vector field on $S$, we denote by $\exp(tX)$ the local
one-parameter group of diffeomorphisms defined by $X$. If $\mathcal{X}(S)$ is
a family of all vector fields on a subcartesian space $S$, the orbit $S_{x}$
of $\mathcal{X}(S)$ through $x$ is given by equation (\ref{accessible}). In
other words,
\[
S_{x}=\{\exp(t_{n}X_{n})\raisebox{2pt}{$\scriptstyle\circ\, $}...\raisebox
{2pt}{$\scriptstyle\circ\, $}\exp(t_{1}X_{1})\mid n\in\mathbb{N},\text{ }%
t_{1},...,t_{n}\in\mathbb{R},\text{ }X_{1},...,X_{n}\in\mathcal{X}(S)\}.
\]

\begin{theorem}
Let $\mathcal{X}(S)$ be the family of all vector fields on a subcartesian
space $S$. For each $x\in S$, the orbit $S_{x}$ is a manifold, and the
inclusion map $S_{x}\hookrightarrow S$ is smooth.
\end{theorem}


\begin{thebibliography}{99}
\bibitem{Andrzejczak}G. Andrzejczak, ``On regular tangent covectors, regular
differential forms, and smooth vector fields on a differential space'',
\textit{Colloq. Math}. \textbf{46} (1982) 243-255.

\bibitem{Aronszajn}N. Aronszajn, ``Subcartesian and subriemannian spaces'',
\textit{Notices Amer. Math. Soc.} \textbf{14} (1967) 111.

\bibitem{aronszajn-szeptycki1}N.\ Aronszajn and P. Szeptycki, ``The theory of
Bessel potentials'', Part IV,\textit{\ Ann.\ Inst. Fourier} (Grenoble)
\textbf{26} (1975) 27-69.

\bibitem{aronszajn-szeptycki}N.\ Aronszajn and P. Szeptycki, ``Subcartesian
spaces'',\textit{\ J. Differential Geom}. \textbf{15} (1980) 393-416.

\bibitem{Bloch-Leonard}A.M. Bloch and N.E. Leonard, ``Symmetries, conservation
laws and control'', in \textit{Geometry, Mechanics and Dynamics, Volume in
Honour of the 60th Birthday of J.E. Marsden}, P. Newton, Ph. Holmes and A.
Weinstein (eds.), Springer Verlag, New York, 2002.

\bibitem{Brockett}R.W. Brockett, ``Control theory and analytical mechanics'',
in \textit{1976 Ames Research Centre (NASA) Conference on Geometric Control
Theory}, R. Hermann and C. Martin (eds.), \textit{Lie Groups: History,
Frontiers and Applications}, \textbf{7,} Math. Sci. Press, Brookline, Mass., USA.

\bibitem{CHMR}H. Cendra, D.D. Holm, J.E. Marsden and T. Ratiu, ``Lagrangian
reduction, the Euler-Poincar\'{e} equation and semidirect products'',
\textit{Trans. Amer. Math. Soc}. \textbf{186} (1998) 1-25.

\bibitem{BC}R. Cushman and L. Bates, \textit{Global aspects of classical
integrable systems}, Birkh\"{a}user, Basel, 1997.

\bibitem{CS}R. Cushman and J. \'{S}niatycki, ``Differential structure of orbit
spaces'', \textit{Canad. J. Math.} \textbf{53} (2001) 715--755.

\bibitem{DK}J.J. Duistermaat and J.A.C. Kolk, \textit{Lie Groups}, Springer
Verlag, New York, 1999.

\bibitem{Grizzle-Marcus}J.W. Grizzle and S. Marcus, ``The structure of
nonlinear control systems possessing symmetries'', \textit{IEEE Trans.
Automatic Control,} \textbf{30} (1985) 248-258.

\bibitem{Kowalczyk}A. Kowalczyk, ``Tangent differential spaces and smooth
forms'', \textit{Demonstratio Math}. \textbf{13} (1980) 893-905.

\bibitem{marshall}C.D.\ Marshall, ``Calculus on subcartesian spaces'',
\textit{J. Differential Geom}. \textbf{10} (1975) 551-573.

\bibitem{marshall1}C.D.\ Marshall, ``The de Rham cohomology on subcartesian
spaces'', \textit{J. Differential Geom}. \textbf{10} (1975) 575-588.

\bibitem{Matuszczyk}H. Matuszczyk, ``On the tangent spaces to differential
spaces'', \textit{Demonstratio Math}. \textbf{14} (1981) 937-942 (1982).

\bibitem{Matuszczyk1}H. Matuszczyk, ``On the tangent and cotangent bundle of a
differential space'', \textit{Ann. Polon. Math.} \textbf{43} (1983) 317-321.

\bibitem{Matuszczyk-Waliszewski}H. Matuszczyk and W. Waliszewski, ``A
nonclassical definition of tangent bundle and cotangent bundle'',
\textit{Demonstratio Math.} \textbf{15} (1982) 913-924.

\bibitem{OR}J.-P. Ortega and T.S. Ratiu, ``The optimal momentum map'', in
\textit{Geometry, Mechanics and Dynamics, Volume in Honour of the 60th
Birthday of J.E. Marsden}, P. Newton, Ph. Holmes and A. Weinstein (eds.),
Springer Verlag, New York, 2002.

\bibitem{Schaft}A.J. van der Schaft, ``Symmetries and conservation laws for
Hamiltonian systems with inputs and outputs: a generalization of Noether's
theorem'', \textit{Syst. Contr. Letters}, \textbf{1} (1981) 108-115.

\bibitem{Schwarz}G.W. Schwarz, ``Smooth functions invariant under the action
of a compact Lie group'', Topology \textbf{14} (1975) 63-68.

\bibitem{sikorski1}R. Sikorski, ``Abstract covariant derivative'',
\textit{Colloq. Math}. \textbf{18} (1967) 251-272.

\bibitem{sikorski2}R. Sikorski, ``Differential modules'', \textit{Colloq.
Math}. \textbf{24} (1971) 45-79.

\bibitem{Sikorski}R.\ Sikorski, \textit{Wst\k{e}p do Geometrii R\'{o}%
\.{z}niczkowej}, PWN, Warszawa, 1972.

\bibitem{SL}R. Sjamaar and E. Lerman, ``Stratified symplectic spaces and
reduction'', \textit{Ann. Math.} \textbf{134} (1991) 375-422.

\bibitem{Sniatycki2001}J. \'{S}niatycki, ``Almost Poisson structures and
nonholonomic singular reduction'', \textit{Rep. Math. Phys}., \textbf{48}
(2001) 235-248.

\bibitem{Sniatycki}J. \'{S}niatycki, ``Orbits of families of vector fields on
subcartesian spaces'', (preprint) arXiv:math.DG/0211212v1.

\bibitem{spallek1}K. Spallek, ``Differenzierbare R\"{a}ume'', \textit{Math.
Ann}., \textbf{180} (1969) 269-296.

\bibitem{spallek2}K. Spallek, Differential forms on differentiable spaces'',
\textit{Rend. Mat.} (2) \textbf{6} (1971) 237-258.

\bibitem{Stefan}P. Stefan: ``Accessible sets, orbits and foliations with
singularities'', \textit{Proc. London Math. Soc}., \textbf{29} (1974) 699-713.

\bibitem{Sussmann}H. J. Sussmann, ``Orbits of families of vector fields and
integrability of distributions'', \textit{Trans. Amer. Math. Soc}.
\textbf{180} (1973) 171-188.

\bibitem{Walczak}P.G. Walczak, ``A theorem on diffeomorphisms in the category
of differential spaces'', \textit{Bull. Acad. Polon. Sci., S\'{e}r. Sci. Astr.
Math. Phys.,} \textbf{21} (1973) 325-329.

\bibitem{Wierzbicki}W. Wierzbicki, ``The tangent and the arcwise tangent space
to a differential space'', \textit{Ann. Polon. Math}. \textbf{40} (1983) 207-212.
\end{thebibliography}
\end{document}